\documentclass[11pt]{article}
\usepackage{amssymb}
\usepackage{graphicx}
\usepackage{amsfonts}
\usepackage{latexsym}
\usepackage{amsthm}
\usepackage{amsmath}

\usepackage{amssymb, amscd}

\usepackage{url}

\newtheorem{problem}{Problem}[section]
\newtheorem{theo}[problem]{Theorem}
\newtheorem{rem}[problem]{Remark}

\newtheorem{defin}[problem]{Definition}
\newtheorem{prop}[problem]{Proposition}
\newtheorem{cor}[problem]{Corollary}

\newtheorem{exam}[problem]{Example}

\begin{document}
\date{}
 \title{{\huge \bf Combinatorial groupoids, \\ cubical complexes, and\\
 the Lov\' asz conjecture}}

\author{Rade  T. \v Zivaljevi\' c}

\maketitle

\begin{abstract}
A foundation is laid for a theory of {\em combinatorial
groupoids}, allowing us to use concepts like ``holonomy'',
``parallel transport'', ``bundles'', ``combinatorial curvature''
etc.\ in the context of simplicial (polyhedral) complexes, posets,
graphs, polytopes and other combinatorial objects. A new,
holonomy-type invariant for cubical complexes is introduced,
leading to a combinatorial ``Theorema Egregium'' for cubical
complexes non-embeddable into cubical lattices. Parallel transport
of $Hom$-complexes and maps is used as a tool for extending
Babson-Kozlov-Lov\'{a}sz graph coloring results to more general
statements about non-degenerate maps (colorings) of simplicial
complexes and graphs.

\end{abstract}

\section{Introduction}

There has been new and encouraging evidence that the language and
methods of the theory of groupoids, after being successfully
applied in other major mathematical fields, offer new insights and
perspectives for applications in combinatorics and discrete and
computational geometry.

The groupoids (groups of projectivities) have recently appeared in
geometric combinatorics in the work of M. Joswig \cite{Josw2001},
see also a related paper with Izmestiev \cite{IzmJos2002} and the
references to these papers, where they have been applied to toric
manifolds, branched coverings over $S^3$, colorings of simple
polytopes etc.

The purpose of this paper is to show that these developments
should not be seen as isolated examples. Quite the opposite, they
serve as a motivation for further extensions and generalizations.
As a first application we show how a cubical extension of Joswig's
groupoid provides new insight about cubical complexes
non-embeddable into cubical lattices (a question related to a
problem of S.P. Novikov which arose in connection with the
$3$-dimensional Ising model) \cite{BuPa02} \cite{Novikov}. The
second application leads to a generalization, from graphs to
simplicial complexes, of a recent resolution of the Lov\' asz
conjecture by Babson and Kozlov \cite{BabsonKozlov2}
\cite{Kozlov-review}.

Combinatorial groupoids also appear (implicitly) in other
contemporary combinatorial constructions  and applications, see
Babson et al.\ \cite{BBLL}, recent work of Barcelo et al. on
Atkin-homotopy \cite{BKLW} \cite{BL} etc. One of our central
objectives is to advocate their systematic use, and  to propose
their recognition as a valuable tool for geometric and algebraic
combinatorics.

\subsection{The unifying theme}

Our point of departure is an observation that different problems,
including the question of existence of an embedding $f : M
\hookrightarrow N$ of Riemann manifolds, the existence and
classification of embeddings (immersions) of cubical surfaces into
hypercubes or cubical lattices, the existence of a coloring of a
graph $G=(V_G,E_G)$ with a prescribed number of colors etc., can
all be approached from a similar point of view.

\begin{figure}[hbt]
\centering
\includegraphics[scale=0.70]{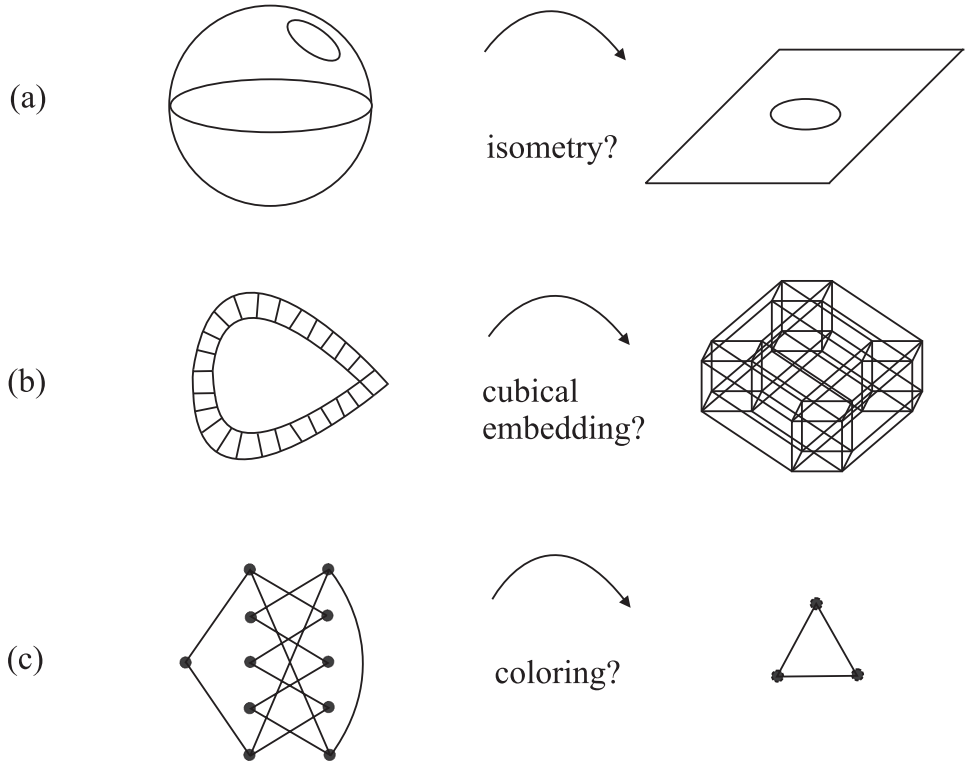}
\caption{What is a common point of view?!} \label{fig:sphere7}
\end{figure}

The unifying theme and a single point of view is provided by the
concept of a groupoid. In each of the listed cases there is a
groupoid naturally associated to an object of the given category.
For each of these groupoids there is an associated ''parallel
transport'', holonomy groups and other related invariants which
serve as obstructions for the existence of  morphisms indicated in
Figure~\ref{fig:sphere7}.

If $M=(M,g)$ is a Riemannian manifold, the associated groupoid
$\mathcal{G}={\mathcal{G}_{M}}$ has $M$ as the set of objects
while the morphism set $\mathcal{G}(x,y)$ consists of all linear
isomorphisms $a : T_xM\rightarrow T_yM$ arising from a parallel
transport along piece-wise smooth curves from $x$ to $y$. One of
the manifestations of Gauss ``Theorema Egregium'' is that the
associated holonomy group is a metric invariant, consequently an
isometry depicted in Figure~\ref{fig:sphere7}~(a)  is not
possible.

In a similar fashion a cubical ``Theorema Egregium''
(Theorem~\ref{thm:cub-thm-egr}) provides an obstruction for an
embedding of a cubical complex into a hypercube or a cubical
lattice. As a consequence the cubulation (quadrangulation) of an
annulus depicted in Figure~\ref{fig:sphere7}~(b) does not admit a
cubical embedding into a hypercube for the same formal reason
(incompatible holonomy) a spherical cap cannot be isometrically
represented in the plane.

Perhaps it comes as a surprise that the graph coloring problem,
Figure~\ref{fig:sphere7}~(c), can be also approached from a
similar point of view. An analysis of holonomies (``parallel
transport'') of diagrams of $Hom$-complexes of graphs (simplicial
complexes) over the associated Joswig groupoid eventually leads to
a general result (Theorem~\ref{thm:main}) which includes the
``odd'' case of Babson-Kozlov-Lov\'{a}sz coloring theorem as a
special case.

\subsection{An overview of the paper}

Our objective is to lay a foundation for a program of associating
groupoids to posets, graphs, complexes, arrangements,
configurations and other combinatorial structures. This
``geometrization of combinatorics'' should provide useful guiding
principles and transparent geometric language introducing concepts
such as combinatorial holonomy, discrete parallel transport,
combinatorial curvature, combinatorial bundles etc. The first
steps in this direction are made in Section~\ref{sec:groupoids}.

Section~\ref{sec:cubical}, devoted to cubical complexes, is an
elaboration of the theme depicted in Figure~\ref{fig:sphere7}~(b).
A holonomy $\mathbb{Z}_2$-valued invariant $I(K)$ of a cubical
complex is introduced. This leads to a ``combinatorial curvature''
$CC(K)$ and associated ``Cubical Theorema Egregium''
(Theorem~\ref{thm:cub-thm-egr}) which provides necessary
conditions for the existence of an embedding/immersion
$K\hookrightarrow L$ of cubical complexes.

Section~\ref{sec:gen-lov-conj}, formally linked to
Figure~\ref{fig:sphere7}~(c), introduces groupoids to the problem
of coloring graphs and complexes. The focus is on the Lov\'{a}sz
$Hom$-conjecture and its ramifications. Parallel transport of
graph and more general $Hom$-complexes over graphs/complexes is
introduced and fundamental invariance of homotopy types of maps is
established in Propositions~\ref{prop:transport} and
\ref{prop:prenos}. This ultimately leads to general results about
coloring of graphs and simplicial complexes
(Theorems~\ref{thm:main} and \ref{thm:paradigm},
Corollary~\ref{cor:LBK-gen}).

\section{Groupoids}\label{sec:groupoids}

It is, or should be, well known that the concept of a {\em group}
is sometimes not sufficient to deal with the concept of symmetry
in general. Groupoids, understood as groups with ``space-like
properties'', often alow us to handle objects which exhibit what
is clearly recognized as symmetry although they admit no global
automorphism whatsoever. Unlike groups, groupoids are capable of
describing reversible processes which can pass through a number of
states. For example according to Connes \cite{Connes}, Heisenberg
discovered quantum mechanics by considering the groupoid of
quantum transitions rather than the group of symmetry.

By definition a groupoid is a small category  $\mathcal{C} =
(Ob(\mathcal{C}), Mor(\mathcal{C}))$ such that each morphism
$\alpha\in Mor(\mathcal{C})$ is an isomorphism. We follow the
usual terminology and notation \cite{Brown87} \cite{Brown88}
\cite{Higg} \cite{Wein96}. Perhaps the only exception is that we
call the {\em vertex} (isotropy) group $\Pi(\mathcal{C},x):=
\mathcal{C}(x,x)$ the {\em holonomy} group of $\mathcal{C}$ at
$x\in Ob(\mathcal{C})$. If the groupoid $\mathcal{C}$ is connected
all holonomy groups $\Pi(\mathcal{C},x)$ are isomorphic and their
isomorphism type is denoted by $\Pi(\mathcal{C})$.

\subsection{Generalities about ``bundles'' and ``parallel transport''}
\label{sec:generalities}

The notion of groupoid is a common generalization of the concepts
of space and group, i.e.\ the theory of groupoids allows us to
treat spaces, groups and objects associated to them from the same
point of view. A common generalization for the concept of a bundle
$Y$ over $X$ and a $G$-space $Y$ is a $\mathcal{C}$-space or more
formally a diagram over the groupoid $\mathcal{C}$ defined as a
functor $F : \mathcal{C}\rightarrow Top$.

In order to preserve the original intuition we use a geometric
language which resembles the usual concepts as much as possible.
In particular we clarify what is in this paper meant by a
$\mathcal{C}$-parallel transport on a ``bundle'' of spaces over a
set $S$.

A (naive) ``bundle'' is a map $\phi : X\rightarrow S$. We assume
that $S$ is a set and that $X(i):=\phi^{-1}(i)$ is a topological
space, so a bundle is just a collection of spaces (fibres) $X(i)$
parameterized by $S$. If all spaces $X(i)$ are homeomorphic to a
fixed ``model'' space, this space is referred to as the fiber of
the bundle $\phi$.

Suppose that $\mathcal{C}$ is a groupoid on $S$ as the set of
objects. In other words $\mathcal{C}=(Ob(\mathcal{C}),
Mor(\mathcal{C}))$ is a small category where $Ob(\mathcal{C})=S$,
such that all morphisms $\alpha\in Mor(\mathcal{C})$ are
invertible.

A ``connection'' or ``parallel transport'' on the bundle
$\mathcal{X}=\{X(i)\}_{i\in S}$ is a functor (diagram)
$\mathcal{P}: \mathcal{C}\rightarrow Top$ such
$X(i)=\mathcal{P}(i)$ for each $i\in S$.

Informally speaking, the groupoid $\mathcal{C}$ provides a ``road
map'' on $S$, while the functor $\mathcal{P}$ defines the
associated transport from one fibre to another.

Sometimes it is convenient to view the bundle
$\mathcal{X}=\{X(i)\}_{i\in S}$ as a map $\mathcal{X}:
S\rightarrow Top$. Then to define a ``connection'' on this bundle
is equivalent to enriching the map $\mathcal{X}$ to a functor
$\mathcal{P}: \mathcal{C}\rightarrow Top$.

\subsection{Combinatorial groupoids}
\label{subsec:comb-group}

The following definition introduces the first in a series of {\em
combinatorial grou\-po\-ids}.

\begin{defin}\label{def:comb-groupod}
Suppose that $(P,\leq)$ is a (not necessarily finite) poset.
Suppose that $\Sigma$ and $\Delta$ two families of subposets of
$P$. Choose $\sigma_1,\sigma_2\in \Sigma$. If for some $\delta\in
\Delta$ both $\delta\subset \sigma_1$ and $\delta\subset
\sigma_2$, then the posets $\sigma_1,\sigma_2$ are called
$\delta$-adjacent, or just adjacent if $\delta$ is not specified.
Define $\mathcal{C}=(Ob(\mathcal{C}), Mor(\mathcal{C}))$ as a
small category over $Ob(\mathcal{C}) = \Sigma$ as the set of
objects as follows. For two $\delta$-adjacent objects $\sigma_1$
and $\sigma_2$, an {\em elementary morphism} $\alpha\in
\mathcal{C}(\sigma_1,\sigma_2)$ is an isomorphism $\alpha :
\sigma_1\rightarrow \sigma_2$ of posets which leaves $\delta$
point-wise fixed. A morphism $\mathfrak{p} \in
\mathcal{C}(\sigma_0,\sigma_m)$ from $\sigma_0$ to $\sigma_m$ is
an isomorphism of posets $\sigma_0$ and $\sigma_m$ which can be
expressed as a composition of elementary morphisms.
\end{defin}

Given two adjacent objects $\sigma_1$ and $\sigma_2$, an
elementary morphism $\alpha\in \mathcal{C}(\sigma_1,\sigma_2)$ may
not exist at all, or if it exists it may not be unique. In case it
exists and is unique it will be frequently denoted by
$\overrightarrow{\sigma_1\sigma_2}$ and sometimes referred to as a
``flip'' from $\sigma_1$ to $\sigma_2$. In this case a morphism
$\mathfrak{p}\in \mathcal{C}(\sigma_0,\sigma_m)$ is by definition
a composition of flips
\[
\mathfrak{p} = \overrightarrow{\sigma_0\sigma_1}\ast
\overrightarrow{\sigma_1\sigma_2}\ast\ldots\ast
\overrightarrow{\sigma_{n-1}\sigma_n}.
\]

\medskip\noindent
{\bf Caveat:} Here we adopt a useful convention that $(x)(f\ast
g)= (g\circ f)(x)$ for each two composable maps $f$ and $g$. The
notation $f\ast g$ is often given priority over the usual $g\circ
f$ if we want to emphasize that the functions act on the points
from the right, that is if the arrows in the associated formulas
point from left to the right.

\medskip

\begin{rem}{\rm
Definition~\ref{def:comb-groupod} can be restated in a much more
general form. For example the ambient poset $P$ can be replaced by
a small category $\mathcal{G}$, $\Sigma$ and $\Delta$ by families
of subobjects of $\mathcal{G}$, elementary morphisms are defined
as commutative diagrams etc. However our objective in this paper
is not to explore all the possibilities. Instead we create an
``ecological niche'' for combinatorial groupoids which may be
populated by new examples and variations as the theory develops.}
\end{rem}

Suppose that $P$ is a ranked poset of depth $n$ with the
associated rank function $r : P\rightarrow [n]$. Let $\mathcal{E}
= \mathcal{E}_P$ be the $\mathcal{C}$-groupoid described in
Definition~\ref{def:comb-groupod} associated to the families
$\Sigma :=\{P_{\leq x} \mid r(x)=n\}$ and $\Delta :=\{P_{\leq y}
\mid r(y)=n-1\}$. It is clear that other ``rank selected''
groupoids can be similarly defined.

\medskip
The definitions of groupoids $\mathcal{C}$ and $\mathcal{E}$ are
easily extended from posets to simplicial, polyhedral, or other
classes of cell complexes. If $K$ is a complex and $P := P_K$ the
associated face poset, then $\mathcal{C}_K$ and
$\mathcal{E}_K=\mathcal{E}(K)$ are groupoids associated to the
poset $P_K$. We will usually drop the subscript whenever it is
clear from the context what is the ambient poset $P$ or complex
$K$.

\begin{exam}{\rm Suppose that $K$ is a pure, $d$-dimensional
simplicial complex. Let $\mathcal{E}(K)$ be the associated
$\mathcal{E}$-groupoid corresponding to ranks $d$ and $d-1$. Then
the {\em groups of projectivities} $\Pi(K,\sigma)$, introduced by
Joswig in \cite{Josw2001}, are nothing but the holonomy groups of
the groupoid $\mathcal{E}(K)$. For this reason the groupoid
$\mathcal{E}(K)$ is in the sequel often referred to as Joswig's
groupoid and denoted by $\mathcal{J}(K)$. $\mathcal{J}(K)$ is {\em
connected} as a groupoid if and only if $K$ is ``strongly
connected'' in the sense of \cite{Josw2001}.
 }
\end{exam}

A simplicial map of simplicial complexes is non-degenerate if it
is 1--1 on simplices. The following definition extends this
concept to the case of posets.

\begin{defin}\label{def:nondeg}
A monotone map of posets $f : P\rightarrow Q$ is {non-degenerate}
if the restriction of $f$ on $P_{\leq x}$ induces an isomorphism
of posets $P_{\leq x}$ and $Q_{\leq f(x)}$ for each element $x\in
P$. Similarly, a map of simplicial, cubical or more general cell
complexes is non-degenerate if the associated map of posets is
non-degenerate. In this case we say that $P$ is {\em mappable} to
$Q$ while a non-degenerate map $f : P\rightarrow Q$ is often
referred to as a {\em combinatorial immersion} from $P$ to $Q$.
\end{defin}

\begin{exam}
{\rm A graph homomorphism \cite{Kozlov-review} $f : G_1
\rightarrow G_2$ can be defined as a non-degenerate map of
associated $1$-dimensional cell complexes. A $n$-coloring of a
graph $G$ is a non-degenerate map (graph homomorphism) $f :
G\rightarrow K_n$ where $K_n$ is a complete graph on $n$
vertices.}
\end{exam}

\begin{prop}\label{prop:functor}
Suppose that $P$ and $Q$ are ranked posets of depth $n$ and let $f
: P\rightarrow Q$ be a non-degenerate map. Then there is an
induced map (functor) $F : \mathcal{E}_P \rightarrow
\mathcal{E}_Q$ of the associated $\mathcal{E}$-groupoids.
Moreover, $F$ induces an inclusion map $\Pi(\mathcal{E}_P,
\mathfrak{p})\hookrightarrow \Pi(\mathcal{E}_Q, F(\mathfrak{p}))$
of the associated holonomy groups.
\end{prop}

\section{Cubical complexes}
\label{sec:cubical}

The following two problems serve as a motivation for studying
non-degenerate morphisms between cubical complexes, the case of
embeddings being of special importance.

\begin{enumerate}

\item[${\mathbf P_1}$] (S.P. Novikov) Characterize $k$-dimensional
complexes that admit a (cubical) embedding, or more generally a
``combinatorial immersion'', into the standard cubical lattice of
$\mathbb{R}^d$ for some $d$, \cite{Novikov} \cite{BuPa02}.

\item[${\mathbf P_2}$] (N. Habegger) Suppose we have two
cubulations of the same manifold. Are they related by the bubble
moves, \cite{Kirby}?
\end{enumerate}

The first question was according to \cite{BuPa02} motivated by a
problem from statistical physics ($3$-dimensional Ising model).
The second initiated the study of ``bubble moves'' on cubical
subdivisions (cubulations) of manifolds and complexes which are
analogs of stellar operations in simplicial category. For these
and related questions about cubical complexes the reader is
referred to \cite{BC} \cite{BEE} \cite{DSS86} \cite{DSS87}
\cite{Epp99} \cite{Fu99} \cite{Fu99b} \cite{Fu05}
\cite{Karalashvili} \cite{SZ04} and the references in these
papers.

\medskip
Recall that a cell complex $K$ is cubical if  it is a regular
$CW$-complex such that the associated face poset $P_K$ is cubical
in the sense of the following definition.

\begin{defin}\label{def:cubical}
$P$ is a cubical poset if:
\begin{enumerate}
 \item[{\rm (a)}] for each $x\in P$, the subposet $P_{\leq x}$ is
isomorphic to the face poset of some cube $I^q;$
 \item[{\rm (b)}] $P$ is a semilattice in the sense that if a pair
 $x,y\in P$ is bounded from above then it has the least upper
 bound.
\end{enumerate}
\end{defin}

If a space $X$ comes equipped with a standard cubulation, clear
from the context, this cubical complex is denoted by $\{X\}$, the
associated $k$-skeleton is denoted by $\{X\}_{(k)}$ etc. For
example $\{I^d\}_{(k)}$ is the $k$-skeleton of the standard
cubulation of the $d$-cube.

The group $BC_k$ of all symmetries of a $k$-cube is isomorphic to
the group of all signed, permutation $(k\times k)$-matrices. Its
subgroup of all matrices with even number of $(-1)$-entries is
denoted by $BC_k^{even}$.

\begin{figure}[hbt]
\centering
\includegraphics[scale=0.50]{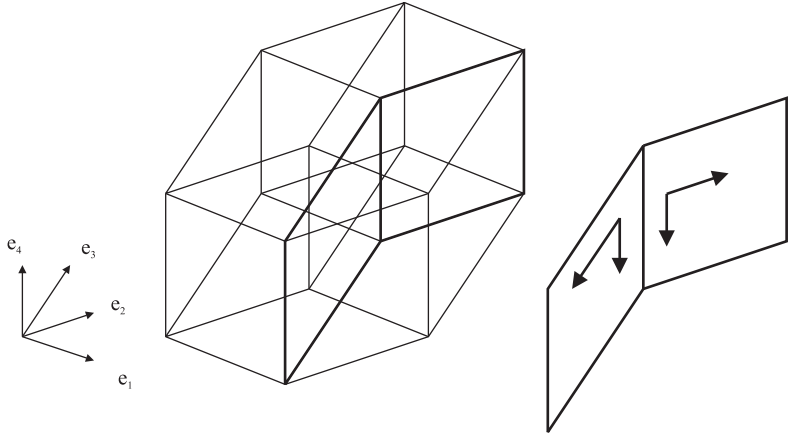}
\caption{The effect of a ``flip'' on the sign characteristic.}
\label{fig:flip}
\end{figure}

\begin{theo}\label{thm:cubical-lattice}
Suppose that $K$ is a $k$-dimensional cubical complex which is
embeddable/mappable to $\{\mathbb{R}^n\}_{(k)}$, the
$k$-dimensional skeleton of the standard cubical decomposition
$\{\mathbb{R}^n\}=\{\mathbb{R}^1\}^n$ of $\mathbb{R}^n$. Then
$\Pi(K,\sigma)\subset BC_k^{even}$ for each cube $\sigma\in K$.
Moreover, if $\sigma\in L\subset K$ where $L\cong
\{I^{k+1}\}_{(k)} $ is a subcomplex of $K$ isomorphic to the
$k$-skeleton of the $(k+1)$-dimensional cube $\{I^{k+1}\}$, then
\[
\Pi(K,\sigma)\cong BC_k^{even}.
\]
\end{theo}

\medskip\noindent
{\bf Proof:} The finiteness of the holonomy group $\Pi(K,\sigma)$
allows us to assume that the complex $K$ is also finite, hence
embeddable/mappable to $\{[0,\nu]^n\}_{(k)}$ for some integer
$\nu\geq 1$. The case $k=1$ of the theorem is a consequence of the
simple fact that each closed edge-path in $\{\mathbb{R}^n\}$ is of
even length, so let us assume that $k\geq 2$.

Each finite cubical subcomplex of $\{\mathbb{R}^n\}_{(k)}$ is
embeddable in a cube $\{I^d\}_{(k)}$ for sufficiently large $d$.
Indeed, the $1$-complex $\{I^d\}_{(1)}$ is a Hamiltonian graph
with $2^d$-vertices which implies that the ``chain''
$\{[a,b]\}\subset \{\mathbb{R}^1\}$ is embeddable in the
$1$-skeleton of the cube $\{I^d\}$ if $b-a\leq 2^d-1$. Since
$\{I^{d_1}\}\times \{I^{d_2}\}\cong \{I^{d_1+d_2}\}$, we observe
that $\{[a,b]^n\}$ is embeddable in $\{I^d\}$ for some $d$, hence
the same holds for their $k$-skeletons.

In light of Proposition~\ref{prop:functor}, the holonomy group
$\Pi(K,\sigma)$ is a subgroup of $\Pi(\{I^d\}_{(k)})$ for some
$d>k$ so it is sufficient to establish the result for the complex
$K=\{I^d\}_{(k)}$ where $d>k\geq 2$.

As a preliminary step let us recall some basic facts about the
$d$-cube $I^d$. Assume that $I^d\subset \mathbb{R}^d$ and let
$[e_1,e_2,\ldots,e_d]$ be the standard orthonormal frame in
$\mathbb{R}^d$. Let $I^d_{x_i=0}=\{\bar{x}\in I^d \mid x_i=0\}$,
respectively $I^d_{x_i=1}=\{\bar{x}\in I^d \mid x_i=1\}$, be the
{\em front}, respectively {\em back} face of $I^d$ in the
direction of the basic vector $e_i$. Given a face
$\sigma\in\{I^d\}$, let
\[
I_0(\sigma):= \{i\in [d] \mid \sigma\subset I^d_{x_i=0}\} \quad
\mbox{\rm and} \quad I_1(\sigma):= \{i\in [d] \mid \sigma\subset
I^d_{x_i=1}\}.
\]
By definition $\sigma = \{\bar{x}\in I^d\mid x_i = 0 \mbox{ {\rm
and} } x_j=1 \mbox{ {\rm for each} } i\in I_0(\sigma) \mbox{ {\rm
and} } j\in I_1(\sigma)\}$. Let $I(\sigma):=
[d]\setminus(I_0(\sigma)\cup I_1(\sigma))$. A cube $\tau\in
\{I^d\}$ is a front face, respectively back face of $\sigma$, if
for some $i\in I(\sigma), \tau = \sigma\cap I^d_{x_i=0}$,
respectively $\tau = \sigma\cap I^d_{x_i=1}$.

Suppose that $\sigma_0$ and $\sigma_1$ are two adjacent
$k$-dimensional faces in $\{I^d\}$ such that $\tau := \sigma_0\cap
\sigma_1$ is their common $(k-1)$-dimensional face. Let
$\overrightarrow{\sigma_0\sigma_1}\in
\mathcal{E}(\sigma_0,\sigma_1)$ be the associated element in the
combinatorial groupoid of $\{I^d\}_{(k)}$,
Section~\ref{subsec:comb-group}. Recall that
$\overrightarrow{\sigma_0\sigma_1} : \sigma_0\rightarrow \sigma_1$
is the isomorphism of cell complexes (posets) which keeps the
subcomplex $\{\tau\}$ fixed. Since $\sigma_0$ and $\sigma_1$ are
faces of a regular cube, there is a unique isometry
$R_{\sigma_0\sigma_1}: \sigma_0\rightarrow \sigma_1$ of faces
$\sigma_0$ and $\sigma_1$ which keeps the common face $\tau$
point-wise fixed. Let $T(\sigma_0)$ and $T(\sigma_1)$ be the
$k$-dimensional linear subspaces in $\mathbb{R}^d$ tangent (that
is parallel) to $\sigma_0$ and $\sigma_1$ respectively, and let
$r_{\sigma_0\sigma_1}: T(\sigma_0)\rightarrow T(\sigma_1)$ be the
corresponding isometry associated to $R_{\sigma_0\sigma_1}$. The
introduction of these maps allows us to pass from the
combinatorial groupoid $\mathcal{E}=\mathcal{E}(\{I^d\}_{(k)})$ to
the isomorphic groupoid
$\mathcal{E}_{iso}=\mathcal{E}_{iso}(\{I^d\}_{(k)})$ of isometries
of $k$-faces in $I^d$.

The following lemma is the essential step in describing the action
of the groupoid $\mathcal{E}_{iso}$ on the set of {\em admissible}
frames on $k$-faces of the cube $I^d$. By definition a frame $f_a
= [a_1,a_2,\ldots ,a_k]$ is admissible for a face $\sigma$ if $a_i
= \epsilon_i e_{\nu_i}$ for each $i=1,\ldots,k$ where
$\epsilon_i\in\{-1,+1\}$ and $(\nu_1,\nu_2,\ldots,\nu_k)$ is some
permutation of the set $I(\sigma)$. The {\em sign characteristic}
$sc(f_a)$ of an admissible frame $f_a$ is the number of negative
signs in the sequence $\epsilon_1,\epsilon_2,\ldots,\epsilon_k$.

\medskip\noindent
{\bf Lemma:} The linear isometry $r_{\sigma_0\sigma_1} :
T(\sigma_0)\rightarrow T(\sigma_1)$ maps a frame $f_a=[a_1,\ldots
a_k]$, admissible for $\sigma_0$, to a frame
$f_b=[b_1,\ldots,b_k]$ admissible for the face $\sigma_1$ where
$b_j=\eta_je_{\mu_j}$ for each $j$, such that $\eta_j\in\{-1,+1\}$
and $(\mu_1,\mu_2,\ldots,\mu_k)$ is a permutation of
$I(\sigma_1)$. Moreover, $sc(f_a)=sc(f_b)$ if $\tau$ is either a
front face for both $\sigma_0$ and $\sigma_1$ or a back face for
both $\sigma_0$ and $\sigma_1$. In the opposite case the parity of
the sign characteristic changes, more precisely
$sc(f_b)=sc(f_a)\pm 1$.

\medskip\noindent
The proof of the lemma is by inspection with Figure~\ref{fig:flip}
illustrating the case when the sign characteristic of a frame is
affected by the flip from $\sigma_0$ to $\sigma_1$.

Suppose that $\sigma_0,\sigma_1,\ldots, \sigma_m=\sigma_0$ is a
sequence of adjacent $k$-faces in $\{I^d\}_{(k)}$ and let
$\mathfrak{p}=\overrightarrow{\sigma_0\sigma_1}\ast
\overrightarrow{\sigma_1\sigma_2}\ast\ldots\ast
\overrightarrow{\sigma_{m-1}\sigma_0}$ be the associated element
in $\Pi(\{I^d\}_{(k)},\sigma_0)$. Let $R: \sigma_0\rightarrow
\sigma_0$ be the isometry associated to $\mathfrak{p}$ and $r :
T(\sigma_0)\rightarrow T(\sigma_0)$ the corresponding linear
isometry.

Assume that $I(\sigma_0)=\{i_1,i_2,\ldots,i_k\}$ which implies
that $[e_{i_1},e_{i_2},\ldots,e_{i_k}]$ is the canonical
admissible frame for $\sigma_0$. By a successive application the
lemma on pairs $(\sigma_j,\sigma_{j+1})$, we observe that
$r[e_{i_1},e_{i_2},\ldots,e_{i_k}]=
[\epsilon_1e_{\nu_1},\epsilon_2e_{\nu_2},\ldots,
\epsilon_ke_{\nu_k}]$ where $(\nu_1,\nu_2,\ldots,\nu_k)$ is a
permutation of the set $I(\sigma_0)$ such that the number of
negative signs in the sequence
$(\epsilon_1,\epsilon_2,\ldots,\epsilon_k)$ must be even. From
here we deduce that $r\in BC_k^{even}$.

One easily proves by inspection that $\Pi(\{I^3\}_{(2)})\cong
BC_3^{even}\cong \mathbb{Z}_2\oplus \mathbb{Z}_2$. This in turn
provides a key step for the proof that
$\Pi(\{I^{k+1}\}_{(k)})\cong BC_k^{even}$ which completes the
proof of the second part of the theorem. $\hfill \square$

\medskip

\begin{figure}[hbt]
\centering
\includegraphics[scale=0.40]{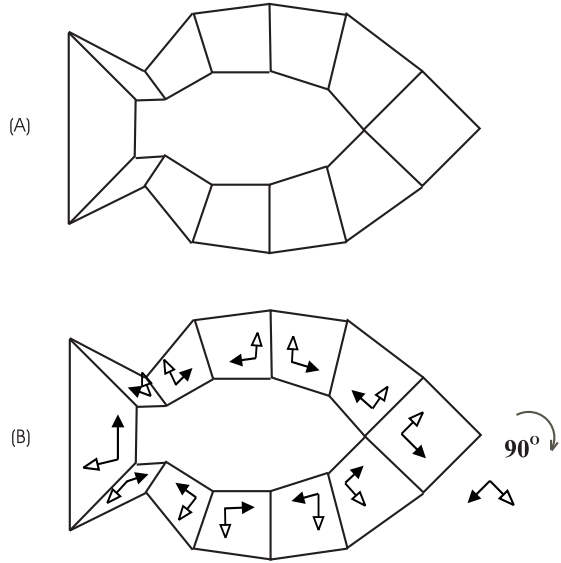}
\caption{Cubical complex non-embeddable into a cubical lattice.}
\label{fig:ribica}
\end{figure}

\begin{cor}\label{cor:ribica}
The complex $K$ depicted in Figure~\ref{fig:ribica} is not
embeddable (mappable) to a cubical lattice (hypercube) of any
dimension. Indeed,
\[
\left[\begin{array}{rl} 0 & 1\\ -1 & 0
\end{array}\right] \in \Pi(K,\sigma)
\]
while by Theorem~\ref{thm:cubical-lattice} only signed permutation
matrices with even number of $(-1)$-entries can arise as
holonomies of subcomplexes of cubical lattices!
\end{cor}

\medskip

Theorem~\ref{thm:cubical-lattice} and Corollary~\ref{cor:ribica}
serve as a motivation for introducing an invariant $I(K)\in
\mathbb{Z}_2$ of a cubical complex $K$ and an associated
''combinatorial curvature'' $CC(K)$ of $K$. Both invariants,
especially the invariant $CC(K)$, can be used for testing if there
exists an embedding (non-degenerate mapping) $f : K\rightarrow L$
of cubical complexes $K$ and $L$. The
$\mathbb{Z}_2$-characteristics $I(K)$ is cruder than $CC(K)$,
however it possesses an additional property that it is invariant
with respect to bubble moves, Proposition~\ref{thm:bubble}.

\begin{defin}
Suppose that $K$ is a $k$-dimensional cubical complex and let
$\Pi(K,\sigma)$ be its combinatorial holonomy group based at
$\sigma\in K$. By definition let $I(K)= 0$ if
$\Pi(K,\sigma)\subset BC_k^{even}$ for all $\sigma$, and $I(K)= 1$
in the opposite case.
\end{defin}

\begin{defin} Given a $k$-dimensional cubical complex, let
Let $Z(K):=\inf\{f_k(W)\mid W\subset K \mbox{ \rm{and} } I(W)=1\}$
where $f_j(K)$ is the number of $j$-di\-men\-sio\-nal cells
(cubes) in $K$. In particular $Z(K)= +\infty$ if $I(K) = 0$. The
{\em combinatorial curvature} of $K$ is the number $CC(K):=
1/Z(K)$ which is by definition $0$ if $I(K) = 0$.
\end{defin}

\begin{rem}{\rm
Both $I(K)$ and $CC(K)$ are global invariants of the cubical
$k$-complex $K$. One can introduce a local invariant $Z(K,\sigma)$
as the minimum length $m$ of a closed chain $\sigma = \sigma_0,
\ldots, \sigma_{m} = \sigma$ of adjacent cubes in $K$ such that
\[
\mathfrak{p}:=\overrightarrow{\sigma_0\sigma_1}\ast
\overrightarrow{\sigma_1\sigma_2}\ast\ldots\ast
\overrightarrow{\sigma_{m-1}\sigma_0}\notin BC_k^{even}.
\]
Then $CC(K,\sigma):=1/Z(K,\sigma)$ is the combinatorial curvature
of $K$ at $\sigma$ and $CC(K)=\sup_{\sigma\in
K^{(k)}}\{C(K,\sigma)\}$. This not only provides a ''correct''
answer in the case of hypercubes (cubical lattices) but resembles
one of the usual definitions of the curvature as the (limit)
quotient of locally defined quantities.}
\end{rem}

The following theorem summarizes the monotonicity properties of
invariants $Z(K)$ and $CC(K)$ in the form suitable for immediate
applications.

\begin{theo}{\rm (Cubical ``Theorema Egregium'')}\label{thm:cub-thm-egr}
If there exists an embedding or more generally a non-degenerate
mapping (cubical immersion) $f : K\rightarrow L$ of
$k$-dimensional cubical complexes $K$ and $L$ then
\[
Z(K)\geq Z(L) \quad \mbox{\rm or equivalently }\quad CC(K)\leq
CC(L).
\]
In other words the ``curvature'' of $K$ must not exceed the
``curvature'' of $L$ if $K$ is to be embedded in $L$.
\end{theo}

\medskip
The following result shows the relevance of the invariant $I(K)$
for the problem of Habegger (problem $P_2$). It says that $I(K)$
is invariant with respect to cubical modifications known as
``bubble moves'', \cite{Fu99} \cite{Fu05}.

\begin{figure}[hbt]
\centering
\includegraphics[scale=0.50]{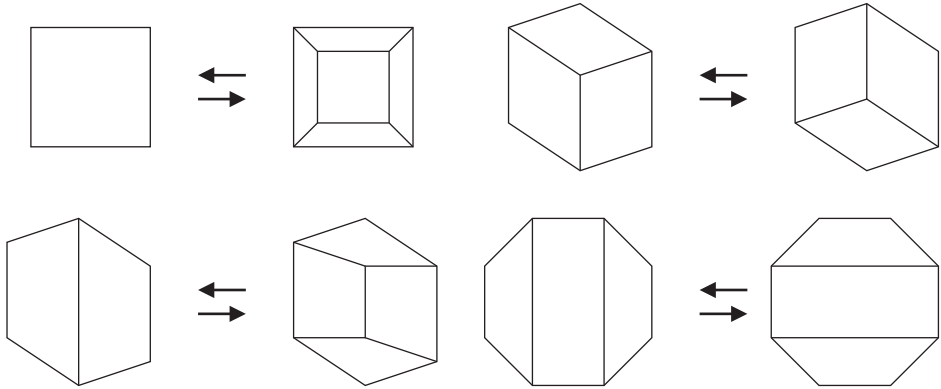}
\caption{Bubble moves for $k=2$.} \label{fig:bubble}
\end{figure}

\begin{prop}\label{thm:bubble} Suppose that $K$ is a $k$-dimensional
cubical complex. Then $I(K)$ is invariant with respect to bubble
moves, that is $I(K) = I(K')$ where $K'$ is the complex obtained
from $K$ by a sequence of bubble moves.
\end{prop}

\medskip\noindent
{\bf Proof:} Suppose that $K'$ is obtained from $K$ by a bubble
move. By definition this means that $K'$ can be obtained from $K$
by excising $B$ from $K$ and replacing it by $B'$ where $B$ and
$B'$ are complementary balls in the boundary of the $(k+1)$-cube,
see Figure~\ref{fig:bubble}. Each cubical path $\mathfrak{p}$ in
$K$ which enters the ball $B\subset K$ can be modified to a
cubical path $\mathfrak{p}'$ by replacing each cubical fragment of
$\mathfrak{p}$ in $B$ by a corresponding fragment in $B'$. Then
the equality $I(K)=I(K')$ follows from $I(\{I^{k+1}\}_{(k)})=0$.
\hfill$\square$

\begin{cor} If $I(K)=1$ and $I(L)=0$ then there does not
exists an embedding (combinatorial immersion) from $K$ to $L$,
even after $K$ and $L$ are modified by some bubble moves.
\end{cor}

\section{Generalized Lov\' asz conjecture}
\label{sec:gen-lov-conj}

In this section we demonstrate how the groupoid of projectivities
introduced by Joswig can be used as a basis for a construction of
parallel transport of graph and more general $Hom$-complexes. In
this framework we develop a general conceptual approach to the
Lov\'{a}sz $Hom$-conjecture, recently resolved by E. Babson and D.
Kozlov \cite{BabsonKozlov2} \cite{Kozlov-review}, and extend their
result both from graphs to simplicial complexes
(Theorem~\ref{thm:main}, Corollary~\ref{cor:LBK-gen}) and to other
test graphs (Theorem~\ref{thm:paradigm}).

\subsection{The Lov\' asz conjecture -- an overview}\label{sec:overview}

One of central themes in topological combinatorics, after the
landmark paper of Laszlo Lov\' asz \cite{Lovasz} where he proved
the classical Kneser conjecture, has been the study and
applications of graph complexes.

The underlying theme is to explore how the topological complexity
of a graph complex $X(G)$ reflects in the combinatorial complexity
of the graph $G$ itself. The results one is usually interested in
come in the form of inequalities $\alpha(X(G))\leq \xi(G)$, or
equivalently in the form of implications
\[
\alpha(X(G))\geq p \Rightarrow \xi(G)\geq q,
\]
where $\alpha(X(G))$ is a topological invariant of $X(G)$, while
$\xi(G)$ is a combinatorial invariant of the graph $G$.

The most interesting candidate for the invariant $\xi$ has been
the chromatic number $\chi(G)$ of $G$, while the role of the
invariant $\alpha$ was played by the ``connectedness'' of $X(G)$,
its equivariant index, the height of an associated characteristic
cohomology class etc., see \cite{Kozlov-review} \cite{Matousek}
\cite{MatZieg} \cite{Ziv2005} for recent accounts.

The famous result of Lov\' asz quoted above is today usually
formulated in the form of an implication
\begin{equation}\label{eq:L78}
Hom(K_2,G) \mbox{ {\rm is $k$-connected} } \Rightarrow \chi(G)\geq
k+3,
 \end{equation}
where $Hom(K_2,G)$ is (together with the neighborhood complex,
``box complex'' etc.) one of avatars of the homotopically unique
$\mathbb{Z}_2$-graph complex of $G$, \cite{MatZieg}
\cite{Ziv2005}. This complex is a special case of a general graph
complex $Hom(H,G)$ (also introduced by L.~Lov\' asz), a cell
complex which functorially depends on the input graphs $H$ and
$G$.

An outstanding conjecture in this area, referred to as the ``Lov\'
asz conjecture'', was that one obtains a better bound if the graph
$K_2$ in (1) is replaced by an odd cycle $C_{2r+1}$. More
precisely Lov\' asz conjectured that
\begin{equation}\label{eq:Lovasz-conj}
Hom(C_{2r+1},G) \mbox{ {\rm is $k$-connected} } \Rightarrow
\chi(G)\geq k+4.
 \end{equation}
This conjecture was confirmed by Babson and Kozlov in
\cite{BabsonKozlov2}, see also \cite{Kozlov-review} for a more
detailed exposition and \cite{Schultz} \cite{Ziv05}\footnote{This
is a preliminary report which served as a basis for the current
Section~\ref{sec:gen-lov-conj}.} for subsequent developments.

Our objective is to develop methods which both offer a simplified
approach to the proof of implication (\ref{eq:Lovasz-conj}), at
least in the case when $k$ is odd, and providing new insight, open
a possibility of proving similar results for other classes of
(hyper)graphs and simplicial complexes.

An example of such a result is Theorem~\ref{thm:main}. One of its
corollaries is the following implication,
\begin{equation}\label{eq:complex-conn}
Hom(\Gamma,K) \mbox{ {\rm is $k$-connected} } \Rightarrow
\chi(K)\geq k+d+3
\end{equation}
which, under a suitable assumption on the complex $\Gamma$ and the
assumption that integer $k$ is odd, extends (\ref{eq:Lovasz-conj})
to the case of pure $d$-dimensional simplicial complexes. A
consequence of (\ref{eq:complex-conn}) is
Theorem~\ref{thm:paradigm} which provides a ``receipt'' how to for
odd $k$ generate new examples of ``homotopy test graphs''
\cite{Kozlov-review}, that is graphs $T$ which satisfy the formula
\begin{equation}\label{eq:test}
Hom(T,G) \mbox{ {\rm is $k$-connected} } \Rightarrow \chi(G)\geq k
+ \chi(T)+ 1.
\end{equation}
Examples include the lower three complexes depicted in
Figure~\ref{fig:piramida} which all have chromatic number $4$ and
clique number $3$.

\subsection{Parallel transport of $Hom$-complexes}

For the definition and basic properties of graph complexes
$Hom(G,H)$ the reader is referred to \cite{Kozlov-review}. More
general $Hom$-complexes associated to simplicial complexes $K$ and
$L$ are introduced in Section~\ref{sec:ramifications}.
Section~\ref{sec:generalities} can be used as a glossary for basic
concepts like ``bundles'', ``parallel transport'', ``connection''
etc.

\subsubsection{Natural bundles and groupoids over simplicial
complexes} \label{sec:natural}

Suppose that $K$ and $L$ are finite simplicial complexes and let
$k$ be an integer such that $0\leq k\leq {\rm dim}(K)$. Denote by
$V(K)=K^{(0)}$ the set of all vertices of $K$. Let $S_k = S_k(K)$
be the set of all $k$-dimensional simplices in $K$. Define a
bundle $\mathcal{F}^L_k : S_k\rightarrow Top$ by the formula
\begin{equation}
\mathcal{F}^L_k(\sigma)=Hom(\sigma, L)\cong Hom(\Delta^{[k+1]},L)
\end{equation}
where $Hom(\sigma, L)$ is one of the $Hom$-complexes introduced in
Section~\ref{sec:Hom} and $\Delta^{[k+1]}$ is the simplex with
$[k+1]=\{1,\ldots, k+1\}$ as the set of vertices. By definition a
typical cell in $Hom(\Delta^{[k+1]},L)$ is of the form
$e=\sigma_1\times\ldots\times\sigma_{k+1}\in L^{k+1}$ where
$\{\sigma_i\}_{i=1}^{k+1}$ is a collection of non-empty simplices
in $L$ such that $\sigma_i\cap\sigma_j=\emptyset$ for $i\neq j$
and
 \[
V(\sigma_1)\ast V(\sigma_2)\ast\ldots\ast V(\sigma_{k+1})\subset
L.
 \]
The corresponding cell in $Hom(\sigma, L)$ is described by a
function $\eta : V(\sigma)\rightarrow
2^{V(L)}\setminus\{\emptyset\}$ such that
$\eta(v_1)\cap\eta(v_2)=\emptyset$ for $v_1\neq v_2$ and
 \[
\ast_{v\in V(\sigma)}\eta(v) = \cup_{v\in V(\sigma)}\eta(v)\in L.
 \]

\medskip\noindent
{\bf Example:} The complex $Hom(\sigma^k,\sigma^m)\cong
Hom(\Delta^{[k+1]},\Delta^{[m+1]})$ is well known in combinatorics
as the $(k+1)$-fold ``deleted product''
$(\Delta^{[m+1]})^{k+1}_\delta$ of $\Delta^{[m+1]}$,
\cite[Section~6.3]{Matousek}. It is well known that for $k=1$ the
associated deleted square $(\Delta^{[m+1]})^{2}_\delta$ is
homeomorphic to a $(m-1)$-dimensional sphere. In other words,
$\mathcal{F}^{\sigma^m}_1 : S_1(K)\rightarrow Top$ is a spherical
bundle naturally associated to the simplicial complex $K$.

\medskip
Our next goal, in the spirit of Section~\ref{sec:generalities}, is
to identify a groupoid on the set $S_k$ which acts on the bundle
$\mathcal{F}^L_k$, i.e.\ a groupoid which provides a parallel
transport of fibres of the bundle $\mathcal{F}^L_k$. It turns out
(Proposition~\ref{prop:connection}) that this groupoid is
precisely the $\mathcal{E}$-groupoid $\mathcal{E}(K^{(k)})$,
associated to the $k$-skeleton of $K$, introduced in
Section~\ref{sec:groupoids}.

It came as a pleasant surprise that this groupoid has already
appeared in geometric combinatorics \cite{Josw2001}
\cite{Jos-Usp}. Indeed, the {\em groups of projectivities}
M.~Joswig introduced and studied in these papers are just the
holonomy groups of a groupoid which we call the $k$-th {\em
groupoid of projectivities} of $K$ and denote by
$\mathcal{J}_k(K)$. In these and subsequent papers \cite{Izm2001}
\cite{Izm-Usp2001} \cite{IzmJos2002}, the groups of projectivities
found interesting applications to toric manifolds, branched
coverings over $S^3$, colorings of simple polytopes, etc. The
following excerpt from \cite{Josw2001} reveals the role this
construction played as a key motivating example for more general
concepts introduced in Section~\ref{sec:groupoids}.
\begin{enumerate}
\item[]
\bigskip\noindent {\small `` ... For each ridge $\rho$ contained in
two facets $\sigma,\tau$, there is a unique vertex
$v(\sigma,\tau)$ which is contained in $\sigma$ but not in $\tau$.
We define the {\em perspectivity} $\langle\sigma, \tau\rangle :
\sigma\rightarrow \tau$ by setting
\[
w\mapsto \left\{ \begin{array}{cl}
                 v(\tau,\sigma) & \mbox{ {\rm if} } w =
                 v(\sigma,\tau)\\
                 w & \mbox{ {\rm otherwise} }
                \end{array}\right.
\]
... The {\em projectivity} $\langle g\rangle$ from $\sigma_1$ to
$\sigma_n$ along $g$ is a concatenation
\[
\langle g\rangle = \langle \sigma_0.\sigma_1, \ldots, \sigma_n
\rangle = \langle\sigma_0, \sigma_1\rangle \langle\sigma_1,
\sigma_2\rangle \ldots \langle\sigma_{n-1}, \sigma_n\rangle
\]
of perspectivities. The map $\langle g\rangle$ is a bijection from
$\sigma_0$ to $\sigma_n$. ... ''}
\end{enumerate}

\smallskip

\begin{prop}\label{prop:connection}
For each simplicial complex $K$ and an auxiliary ``coefficient''
complex $L$, there exists a canonical
$\mathcal{J}_k(K)$-connection $\mathcal{P}^L =
\mathcal{P}^L_{K,k}$ on the bundle $\mathcal{F}^L_k$. In other
words the function $\mathcal{F}^L_k: S_k\rightarrow Top$ can be
enriched (extended) to a functor
 \[\mathcal{F}^L_k: \mathcal{J}_k(K)\rightarrow Top\]
where $\mathcal{J}_k(K)\cong \mathcal{E}_k(K)$ is the $k$-th
Joswig's groupoid of projectivities of $K$.
\end{prop}

\medskip\noindent
{\bf Proof:} If $\overrightarrow{\sigma_0\sigma_1}$ is a
perspectivity from $\sigma_0$ to $\sigma_1$ (a ``flip'' in the
language of Section~\ref{subsec:comb-group}) and if $\eta :
V(\sigma_1) \rightarrow 2^{V(L)}\setminus\{\emptyset\}$ is a cell
in $Hom(\sigma,L)$, then $\mathcal{P}^L : \mathcal{F}^L(\sigma_1)
\rightarrow \mathcal{F}^L(\sigma_0)$ is the map defined by
$\mathcal{P}^L(\overrightarrow{\sigma_0\sigma_1})(\eta):=
\overrightarrow{\sigma_0\sigma_1}\ast\eta$. More generally, if
$\mathfrak{p} = \overrightarrow{\sigma_0\sigma_1}\ast
\overrightarrow{\sigma_1\sigma_2}\ast\ldots\ast
\overrightarrow{\sigma_{n-1}\sigma_n}$ is a ``projectivity''
between $\sigma_0$ and $\sigma_n$, then
\begin{equation}
\mathcal{P}^L(\mathfrak{p}) =
\mathcal{P}^L(\overrightarrow{\sigma_0\sigma_1})\ast
\mathcal{P}^L(\overrightarrow{\sigma_1\sigma_2})\ast\ldots\ast
\mathcal{P}^L(\overrightarrow{\sigma_{n-1}\sigma_n})
\end{equation}
or in other words
\begin{equation}
\mathcal{P}^L(\mathfrak{p})(\eta) =
\overrightarrow{\sigma_0\sigma_1}\ast
\overrightarrow{\sigma_1\sigma_2}\ast\ldots\ast
\overrightarrow{\sigma_{n-1}\sigma_n}\ast\eta =
\mathfrak{p}\ast\eta
\end{equation}
which demonstrates in passing that the map
$\mathcal{P}^L(\mathfrak{p})$ depends only on the morphism
$\mathfrak{p}: \sigma_0\rightarrow \sigma_n$ alone, and not on the
associated path $\sigma_0,\ldots \sigma_n$. \hfill$\square$

\subsubsection{Parallel transport of graph complexes}
\label{sec:functor}

The main motivation for introducing the parallel transport of
$Hom$-complexes is the Lov\'{a}sz conjecture and its
ramifications. This is the reason why the case of graphs and the
graph complexes deserves a special attention. Additional
justification for emphasizing graphs comes from the fact that
graph complexes $Hom(G,H)$ have been studied in numerous papers
and today form a well established part of graph theory and
topological combinatorics. The situation with simplicial complexes
is quite the opposite. In order to extend the theory of
$Hom$-complexes from graphs to the category of simplicial
complexes, many concepts should be generalized and the
corresponding facts established in a more general setting. One is
supposed to recognize the main driving forces and to isolate the
most desirable features of the theory. A result should be a
dictionary/glossary of associated concepts, cf.\ Table~1.
Consequently, Section~\ref{sec:functor} should be viewed as an
important preliminary step, leading to the more general theory
developed in Sections~\ref{sec:ramifications} and \ref{sec:main}.

In order to simplify the exposition we assume, without a serious
loss of generality, that all graphs $G=(V(G),E(G))$ are without
loops and multiple edges. In short, graphs are $1$-dimensional
simplicial complexes. Let $G_{\overline{xy}}\cong K_2$ be the
restriction of $G$ on the edge $\overline{xy}\in E(G)$.

\medskip
Following the definitions from Section~\ref{sec:natural} the map
\[
\mathcal{F}^H : E(G) \longrightarrow Top,
\]
where $\mathcal{F}^H(\overline{xy})= \mathcal{F}^H_{\overline{xy}}
:= Hom(G_{\overline{xy}}, H)$, can be thought of as a ``bundle''
over the graph $G$, with $\mathcal{F}^H_{\overline{xy}}$ in the
role of the ``fibre'' over the edge $\overline{xy}$. More
generally, given a class $\mathcal{C}$ of subgraphs of $G$, say
the subtrees, the chains, the $k$-cliques etc., one can define an
associated ``bundle'' $\mathcal{F}^H_{\mathcal{C}}:
\mathcal{C}\rightarrow Top$ by a similar formula
$\mathcal{F}^H_{\mathcal{C}}(\Gamma):= Hom(\Gamma, H)$, where
$\Gamma\in \mathcal{C}$.

The parallel transport $\mathcal{P}^H$, for a given graph
($1$-dimensional, simplicial complex) $H$, is a specialization of
the parallel transport $\mathcal{P}^L$ introduced in
Section~\ref{sec:natural}. For example if
$\overrightarrow{e_1e_2}$ is the flip (perspectivity) between
adjacent edges $e_1=\overline{x_0x_1}$ and $e_2=\overline{x_1x_2}$
in $G$, and if $\eta : \{x_1,x_2\}\rightarrow
2^{V(H)}\setminus\{\emptyset\}$ is a cell in
$\mathcal{F}^H_{\overline{x_1x_2}} = Hom(G_{\overline{x_1x_2}},
H)$, then $\eta':= \mathcal{P}^H(\overrightarrow{e_1e_2})(\eta) :
\{x_0,x_1\}\rightarrow 2^{V(H)}\setminus\{\emptyset\}$ is defined
by
\[
\eta'(x_0) := \eta(x_2)  \mbox{ {\rm and} } \eta'(x_1) :=
\eta(x_1).
\]

\medskip\noindent
{\bf Fundamental observation:} The constructions of the
connections $\mathcal{P}^L$, respectively $\mathcal{P}^H$, are
quite natural and elementary but it is
Proposition~\ref{prop:transport}, respectively its more general
relative Proposition~\ref{prop:prenos}, that serve as an actual
justification for the introduction of these objects.
Proposition~\ref{prop:transport} allows us to analyze the parallel
transport of homotopy types of maps from the complex $Hom(G,H)$ to
complexes $Hom(G_e,H)$, where $e\in E(G)$, providing a key for a
resolution of the Lov\'{a}sz conjecture in the case when $k$ is an
odd integer.

\medskip\noindent
Implicit in the proof of Proposition~\ref{prop:transport} is the
theory of {\em folds} of graphs and the analysis of natural
morphisms between graph complexes $Hom(T,H)$, where $T$ is a tree,
as developed in \cite{BabsonKozlov1} \cite{Kozl2004}
\cite{Kozlov99}. This theory is one of essential ingredients in
the Babson and Kozlov spectral sequence approach to the solution
of Lov\'{a}sz conjecture. Some of these results are summarized in
Proposition~\ref{prop:help}, in the form suitable for application
to Proposition~\ref{prop:transport}.

As usual $L_m$ is the graph-chain of vertex-length $m$, while
$L_{x_1\ldots x_m}$ is the graph isomorphic to $L_m$ defined on a
linearly ordered set of vertices $x_1,\ldots,x_m$. In this context
the ``chain-flip'' is a generic name for the automorphism $\sigma
: L_{x_1\ldots x_m}\rightarrow L_{x_1\ldots x_m}$ of the
graph-chain such that $\sigma(x_j)=x_{m-j+1}$ for each $j$.

\begin{prop}\label{prop:help}
Suppose that $e_1=\overline{x_0x_1}$ and $e_2=\overline{x_1x_2}$
are two distinct, adjacent edges in the graph $G$. Let $\sigma :
L_{x_0x_1x_2}\rightarrow L_{x_0x_1x_2}$ be the chain-flip
automorphism of $L_{x_0x_1x_2}$ and $\widehat{\sigma}$ the
associated auto-homeomorphism of $Hom(L_{x_0x_1x_2},H)$. Suppose
that $\gamma_{ij}: L_{x_ix_j}\rightarrow L_{x_0x_1x_2}$ is an
obvious embedding and $\widehat{\gamma}_{ij}$ the associated maps
of graph complexes. Then,
 \begin{enumerate}
 \item[{\rm (a)}] the induced map $\widehat{\sigma}:
Hom(L_{x_0x_1x_2},H)\rightarrow Hom(L_{x_0x_1x_2},H)$ is homotopic
to the identity map $I$, and
 \item[{\rm (b)}]  the diagram
\[
\begin{CD}
Hom(L_{x_0x_1x_2},H) @>=>> Hom(L_{x_0x_1x_2},H) \\
 @V\widehat{\gamma}_{01}VV @VV\widehat{\gamma}_{12}V \\
 Hom(L_{x_0x_1},H) @<<\mathcal{P}^H(\overrightarrow{e_1e_2})<
 Hom(L_{x_1x_2},H)
 \end{CD}
\]
 is commutative up to homotopy.
 \end{enumerate}
\end{prop}

\medskip\noindent
{\bf Proof:} Both statements are corollaries of Babson and Kozlov
analysis of complexes $Hom(T,H)$, where $T$ is a tree, and
morphisms $\widehat{e}: Hom(T,H)\rightarrow Hom(T',H)$, where $T'$
is a subtree of $T$ and $e: T'\rightarrow T$ the associated
embedding.

Our starting point is an observation that both $L_{x_0x_1}$ and
$L_{x_1x_2}$ are retracts of the graph $L_{x_0x_1x_2}$ in the
category of graphs and graph homomorphisms. The retraction
homomorphisms $\phi_{ij}: L_{x_0x_1x_2}\rightarrow L_{x_ix_j}$,
where $\phi_{01}(x_0)=x_0, \phi_{01}(x_1)=x_1, \phi_{01}(x_2)=x_0$
and $\phi_{12}(x_0)=x_2, \phi_{12}(x_1)=x_1, \phi_{12}(x_2)=x_2$
are examples of {\em foldings} of graphs. By the general theory
\cite{BabsonKozlov1} \cite{Kozl2004}, the maps
$\widehat{\gamma}_{ij}: Hom(L_{x_0x_1x_2},H)\rightarrow
Hom(L_{x_ix_j},H)$ and $\widehat{\phi}_{ij}:
Hom(L_{x_ix_j},H)\rightarrow Hom(L_{x_0x_1x_2},H)$ are homotopy
equivalences. Actually $\widehat{\gamma}_{ij}$ is a deformation
retraction and $\widehat{\phi}_{ij}$ is the associated embedding
such that $\widehat{\gamma}_{ij}\circ \widehat{\phi}_{ij}=I$ is
the identity map.

The part (a) of the proposition is an immediate consequence of the
fact that $\phi_{01}\circ\sigma\circ\gamma_{01} :
L_{x_0x_1}\rightarrow L_{x_0x_1}$ is an identity map. It follows
that
$\widehat{\gamma_{01}}\circ\widehat{\sigma}\circ\widehat{\phi}_{01}=I$,
and in light of the fact that $\widehat{\gamma_{01}}$ and
$\widehat{\phi}_{01}$ are homotopy inverses to each other, we
conclude that $\widehat{\sigma}\simeq I$.

For the part (b) we begin by an observation that $\phi_{12}\circ
\sigma\circ \gamma_{01} = \overrightarrow{e_1e_2}$. Then,
$\mathcal{P}^H(\overrightarrow{e_1e_2}) =
\widehat{\gamma}_{01}\circ\widehat{\sigma}\circ\widehat{\phi}_{12}$,
and as a consequence of $\widehat{\sigma} \simeq I$ and the fact
that $\widehat{\phi}_{12}\circ \widehat{\gamma}_{12}\simeq I$, we
conclude that
\[
\mathcal{P}^H(\overrightarrow{e_1e_2})\circ \widehat{\gamma}_{12}
= \widehat{\gamma}_{01}\circ\widehat{\sigma}\circ
\widehat{\phi}_{12}\circ\widehat{\gamma}_{12}\simeq
\widehat{\gamma}_{01}.
\]
\hfill $\square$

\begin{prop}\label{prop:transport} Suppose that $x_0,x_1,x_2$ are
distinct vertices in $G$ such that $\overline{x_0x_1},
\overline{x_1x_2}\in E(G)$. Let $\alpha_{ij}:
G_{x_ix_{j}}\rightarrow G$ be the inclusion map of graphs and
$\widehat{\alpha}_{ij}$ the associated map of $Hom(\,\cdot\,, H)$
complexes. Then the following diagram commutes up to a homotopy,

\begin{equation}\label{eqn:diagram1}
\begin{CD}
Hom(G,H) @>=>> Hom(G,H) \\
 @V\widehat{\alpha}_{01}VV @VV\widehat{\alpha}_{12}V \\
 Hom(G_{x_0x_1},H) @<<\mathcal{P}^H(\overrightarrow{e_1e_2})<
 Hom(G_{x_1x_2},H)
 \end{CD}
\end{equation}
\end{prop}

\medskip\noindent
{\bf Proof:} The diagram (\ref{eqn:diagram1}) can be factored as
\begin{equation}\label{eqn:diagram2}
\begin{CD}
  Hom(G,H) @>=>> Hom(G,H) \\
  @V\widehat{\beta} VV @VV\widehat{\beta} V \\
  Hom(G_{x_0x_1x_2},H) @>=>> Hom(G_{x_0x_1x_2},H)\\
  @V\widehat{\gamma}_{01}VV @VV\widehat{\gamma}_{12}V\\
  Hom(G_{x_0x_1},H) @<<\mathcal{P}^H(\overrightarrow{e_1e_2})<
  Hom(G_{x_1x_2},H)
 \end{CD}
\end{equation}
where $\beta$ and $\gamma_{ij}$ are obvious inclusions of
indicated graphs such that $\alpha_{ij}=\beta\circ\gamma_{ij}$.
Then the result is a direct consequence of
Proposition~\ref{prop:help}, part (b).
 \hfill $\square$

\subsection{Babson-Kozlov-Lov\'{a}sz result for odd $k$}

The proof \cite{BabsonKozlov2} of Lov\' asz conjecture splits into
two main branches, corresponding to the parity of a parameter $n$,
where $n$ is an integer which enters the stage as the size of the
vertex set of the complete graph $K_n$.

The first branch relies on Theorem~2.3.\ (loc.\ cit.), more
precisely on part (b) of this result, while the second branch is
founded on Theorem~2.6. Both theorems are about the topology of
the graph complex $Hom(C_{2r+1}, K_n)$. Theorem~2.3.~(b) is a
statement about the height of the first Stiefel-Whitney
class\footnote{Subsequently C.~Schultz discovered \cite{Schultz} a
powerful new way of evaluating this invariant, leading to a much
shorter proof of Lov\'{a}sz conjecture.}, equivalently the
Conner-Floyd index \cite{ConnerFloyd} of the $\mathbb{Z}_2$-space
$Hom(C_{2r+1}, K_n)$. Theorem~2.6. claims that for $n$ even,
$2\iota^\ast_{K_n}$ is a zero homomorphism where
\begin{equation}\label{eqn:thm-2.6.}
 \iota^\ast_{K_n} : \widetilde{H}^\ast(Hom(K_2,K_n);\mathbb{Z})
 \longrightarrow \widetilde{H}^\ast(Hom(C_{2r+1},K_n);\mathbb{Z})
\end{equation}
is the homomorphism  associated to the continuous map
\[\iota_{K_n} : Hom(C_{2r+1},K_n)\rightarrow Hom(K_2,K_n),\] which
in turn comes from the inclusion $K_2\hookrightarrow C_{2r+1}$.

\medskip
The central idea of our approach is an observation that
Theorem~2.6.\ can be incorporated into a more general scheme,
involving the ``parallel transport'' of graph complexes over
graphs.

\begin{theo}\label{thm:glavna}
Suppose that $\alpha : K_2 \rightarrow C_{2r+1}$ is an inclusion
map, $\beta : K_2 \rightarrow K_2$ a nontrivial automorphism of
$K_2$, and
 \[ \widehat{\alpha} : Hom(C_{2r+1},H)\rightarrow
Hom(K_2,H), \, \widehat{\beta} : Hom(K_2,H)\rightarrow Hom(K_2,H)
 \]
the associated maps of graph complexes. Then the following diagram
is commutative up to a homotopy
\begin{equation}\label{eqn:diagram3}
\begin{CD}
  Hom(C_{2r+1},H) @>=>> Hom(C_{2r+1},H) \\
  @V\widehat{\alpha} VV @VV\widehat{\alpha} V \\
  Hom(K_2,H) @<\widehat{\beta} << Hom(K_2,H)
   \end{CD}
\end{equation}
\end{theo}

\medskip\noindent
{\bf Proof:} Assume that the consecutive vertices of $G=C_{2r+1}$
are $x_0, x_1, \ldots , x_{2r}$ and let $e_i =
\overline{x_{i-1}x_i}$ be the associated sequence of edges where
by convention $e_{2r+1}=\overline{x_{2r}x_0}$. Identify the graph
$K_2$ to the subgraph $G_{x_0x_1}$ of $G=C_{2r+1}$.

By iterating Proposition~\ref{prop:transport} we observe that the
diagram
\begin{equation}\label{eqn:diagram4}
\begin{CD}
Hom(C_{2r+1},H) @>=>> Hom(C_{2r+1},H) \\
 @V\widehat{\alpha} VV @VV\widehat{\alpha} V \\
 Hom(G_{x_0x_1},H) @<<\mathcal{P}^H(\mathfrak{p})<
 Hom(G_{x_0x_1},H)
 \end{CD}
\end{equation}
is commutative up to a homotopy, where $\mathfrak{p}=
\overrightarrow{e_1e_2}\ast\ldots\ast\overrightarrow{e_{2r+1}e_1}$.
The proof is completed by the observation that $\mathfrak{p}=
\beta$ in the groupoid $\mathcal{J}(C_{2r+1})\cong
\mathcal{{E}}(C_{2r+1})$. \hfill $\square$

\medskip
Theorem~2.6. from \cite{BabsonKozlov2}, the key for the proof of
Lov\'{a}sz conjecture for odd $k$, is an immediate consequence of
Theorem~\ref{thm:glavna}.

\begin{cor}{\rm (\cite{BabsonKozlov2}, T.2.6.)}
If $n$ is even then $2\cdot\iota^\ast_{K_n}$ is a $0$-map where
$\iota^\ast_{K_n}$ is the map described in line\ {\rm
(\ref{eqn:thm-2.6.})}.
\end{cor}

\medskip\noindent
{\bf Proof:} It is sufficient to observe that for $H=K_n$,
$Hom(K_2,K_n)\cong S^{n-2}$ is an even dimensional sphere and that
the automorphism $\widehat{\beta}$ from the diagram
(\ref{eqn:diagram3}) is in this case essentially an antipodal map.
It follows that $\widehat{\beta}$ changes the orientation of
$Hom(K_2,K_n)$ and as a consequence $\iota^\ast_{K_n}=
-\iota^\ast_{K_n}$. \hfill $\square$

\subsection{Generalizations and ramifications}
\label{sec:ramifications}

In this section we extend the results from
Section~\ref{sec:functor} to the case of simplicial complexes.
This generalization is based on the following basic principles.

Graphs are viewed as $1$-dimensional simplicial complexes. Graph
homomorphisms are special cases of {\em non-degenerate} simplicial
maps of simplicial complexes, see \cite{IzmJos2002}
\cite{Josw2001} and Definition~\ref{def:nondeg}. The definition of
$Hom(G,H)$ is extended to the case of $Hom$-complexes $Hom(K,L)$
of simplicial complexes $K$ and $L$. The groupoids needed for the
definition of the parallel transport of $Hom$-complexes are
described in Section~\ref{sec:natural}. Theory of folds for graph
complexes \cite{BabsonKozlov1} \cite{Kozl2004} is extended in
Section~\ref{sec:tree-like} to the case of $Hom$-complexes in
sufficient generality to allow ``parallel transport'' of homotopy
types of maps between graph complexes. This development eventually
leads to Theorem~\ref{thm:main} which extends
Theorem~\ref{thm:glavna} to the case of $Hom$-complexes $Hom(K,L)$
and represents the currently final stage in the evolution of
Theorem~2.6. from \cite{BabsonKozlov2}.

\begin{table}[hbt]\label{tab:tabela}
\begin{center}
\begin{tabular}{||l|l||}\hline
\multicolumn{2}{||c||}{\bf Dictionary}\\ \hline\hline
 graphs &  simplicial complexes \\ \hline
 trees & tree-like complexes \\ \hline
 foldings of graphs & vertex collapsing of complexes \\ \hline
 graph homomorphisms & non-degenerate simplicial maps \\ \hline
 $Hom(G,H)$ & $Hom(K,L)$ \\ \hline
 chromatic number $\chi(G)$ & chromatic number $\chi(K)$ \\ \hline
\end{tabular}
\end{center}
\caption{Graphs vs.\ simplicial complexes.}
\end{table}

\subsubsection{From $Hom(G,H)$ to $Hom(K,L)$} \label{sec:Hom}

Suppose that $K\subset 2^{V(K)}$ and $L\subset 2^{V(L)}$ are two
(finite) simplicial complexes, on the sets of vertices $V(K)$ and
$V(L)$ respectively.

\begin{defin} The set of all non-degenerate simplicial maps from
$K$ to $L$ is denoted by $Hom_0(K,L)$ where a map $f: K\rightarrow
L$ is {\em non-degenerate} (Definition~\ref{def:nondeg}) if it is
injective on simplices.
\end{defin}

\begin{defin}\label{def:hom(k,l)}
$Hom(K,L)$ is a cell complex with the cells indexed by the
functions $\eta : V(K)\rightarrow 2^{V(L)}\setminus\{\emptyset\}$
such that
 \begin{enumerate}
 \item[{\rm (1)}] for each two vertices $u\neq v$, if $\{u,v\}\in K$ then
$\eta(u)\cap\eta(v)=\emptyset$,
 \item[{\rm (2)}] for each simplex $\sigma\in K$, the join
 ${{\ast}}_{v\in V(\sigma)}~\eta(v)\subset \Delta^{V(L)}$
 of all sets (or $0$-dimensional complexes) $\eta(v)$ is a subcomplex of $L$.
\end{enumerate}
More precisely, each function $\eta$ satisfying conditions {\rm
(1)} and {\rm (2)} defines a cell $c_\eta:= \prod_{v\in
V(K)}~\Delta^{\eta(v)}$ in $Hom(K,L)\subset\prod_{v\in
V(K)}\Delta^{V(L)}$ where by definition $\Delta^S$ is an
(abstract) simplex spanned by vertices in $S$.
\end{defin}

The following example shows that $Hom(K,L)$ are close companions
of graph complexes $Hom(G,H)$.

\begin{exam}\label{exam:companion} {\rm The definition of the complex
$Hom(K,L)$ is a natural extension of $Hom(G,H)$ and reduces to it
if $K$ and $L$ are $1$-dimensional complexes. Moreover,
\[
Hom(G,H)\cong Hom(Clique(G),Clique(H))
\]
where $Clique(\Gamma)$ is the simplicial complex of all cliques in
a graph $\Gamma$. }
\end{exam}

\begin{rem}{\rm
The set $Hom_0(K,L)$ is easily identified as the $0$-dimensional
skeleton of the cell-complex $Hom(K,L)$. Moreover, the reader
familiar with \cite{Kozlov-review} can easily check that
$Hom(K,L)$ is determined by the family $M=Hom_0(K,L)$ in the sense
of Definition~2.2.1. from that paper. }
\end{rem}

\subsubsection{Functoriality of $Hom(K,L)$} \label{sec:naturality}

The construction of $Hom(K,L)$ is functorial in the sense that if
$f : K\rightarrow K'$ is a non-degenerate simplicial map of
complexes $K$ and $K'$, then there is an associated continuous map
$\widehat{f} : Hom(K',L)\rightarrow Hom(K,L)$ of $Hom$-complexes.
Indeed, if $\eta : V(K')\rightarrow
2^{V(L)}\setminus\{\emptyset\}$ is a multi-valued function
indexing a cell in $Hom(K',L)$, then it is not difficult to check
that $\eta\circ f : V(K)\rightarrow
2^{V(L)}\setminus\{\emptyset\}$ is a cell in $Hom(K,L)$.

Even more important is the functoriality of $Hom(K,L)$ with
respect to the second variable since this implies the
functoriality of the bundle $\mathcal{F}^L_k$.

\begin{prop}\label{prop:functoriality}
Suppose that $g : L \rightarrow L'$ is a non-degenerate,
simplicial map of simplicial complexes $L$ and $L'$. Then there
exists an associated map
\[
\widehat{g} : Hom(K,L)\rightarrow Hom(K,L').
\]
\end{prop}

\medskip\noindent
{\bf Proof:} Assume that $\eta : V(K)\rightarrow
2^{V(L)}\setminus\{\emptyset\}$ is a cell in $Hom(K,L)$. Then
$g\circ\eta : V(K)\rightarrow 2^{V(L')\setminus\{\emptyset\}}$ is
a cell in $Hom(K,L')$. Suppose $u$ and $v$ are distinct vertices
in $V(K)$. By assumption $\eta(u)\cap \eta(v)=\emptyset$. We
deduce from here that $g(\eta(u))\cap g(\eta(v))\neq \emptyset$,
otherwise $g$ would be a degenerated simplicial map. The second
condition from Definition~\ref{def:hom(k,l)} is checked by a
similar argument. \hfill $\square$

\subsubsection{Chromatic number $\chi(K)$ and its relatives}

The chromatic number $\chi(K)$ of a simplicial complex $K$ is
\[
\inf\{m\in \mathbb{N}\mid Hom_0(K,\Delta^{[m]})\neq\emptyset\}.
\]
In other words $\chi(K)$ is the minimum number $m$ such that there
exists a non-degenerate simplicial map $f : K\rightarrow
\Delta^{[m]}$. It is not difficult to check that
$\chi(K)=\chi(G_K)$ where $G_K = (K^{(0)}, K^{(1)})$ is the
vertex-edge graph of the complex $K$. In particular $\chi(K)$
reduces to the usual chromatic number if $K$ is a graph, that is
if $K$ is a $1$-dimensional simplicial complex.

\medskip
Aside from the usual chromatic  number $\chi(G)$, there are many
related colorful graph invariants \cite{GodRoy}
\cite{Kozlov-review}. Among the best known are the fractional
chromatic number $\chi_f(G)$ and the circular chromatic number
$\chi_c(G)$ of $G$. These and other related invariants are
conveniently defined in terms of graph homomorphisms into graphs
chosen from a suitable family $\mathcal{F} =\{G_i\}_{i\in I}$ of
test graphs. Motivated by this we offer an extension of the
chromatic number $\chi(K)$ in hope that some genuine, new
invariants of simplicial complexes may arise this way.

\begin{defin}
Suppose that $\mathcal{F}=\{T_i\mid i\in I\}$ is a family of
``test'' simplicial complexes and let $\phi : I\rightarrow
\mathbb{R}$ is a real-valued function. A $T_i$-coloring of $K$ is
just a non-degenerate simplicial map $f:K\rightarrow T_i$ and
$\chi_{(\mathcal{F},\phi)}(K)$, the $(\mathcal{F},\phi)$-chromatic
number
 of $K$, is defined as the infimum
of all weights $\phi(i)$ over all $T_i$-colorings,

\[
\chi_{(\mathcal{F},\phi)}(K):= \inf\{\phi(i) \mid
Hom_0(K,T_i)\neq\emptyset\}.
\]

\end{defin}

\subsubsection{Tree-like simplicial complexes}
\label{sec:tree-like}

The tree-like or vertex collapsible complexes are intended to play
in the theory of $Hom(K,L)$-complexes the role similar to the role
of trees in the theory of graph complexes $Hom(G,H)$.

A pure, $d$-dimensional simplicial complex $K$ is {\em shellable}
\cite{Bjo-et-al} \cite{Zieg-book}, if there is a linear order
$F_1, F_2, \ldots , F_m$ on the set of its facets, such that for
each $j\geq 2$, the complex $F_j\cap(\bigcup_{i<j}~F_i)$ is a pure
$(d-1)$-dimensional subcomplex of the simplex $F_j$. The {\em
restriction} $R_j$ of the facet $F_j$ is the minimal new face
added to the complex $\bigcup_{k<j}~F_k$ by the addition of the
facet $F_j$. Let $r_j := {\rm dim}(R_j)\in\{0,1,\ldots,d\}$ be the
type of the facet $F_j$. If $r_j\neq d$ for each $j$ then the
complex $K$ is collapsible. The collapsing process is just the
shelling order read in the opposite direction. From this point of
view, $R_j$ can be described as a free face in the complex
$\bigcup_{i\leq j}~F_i$, and the process of removing all faces $F$
such that $R_j\subset F\subset F_j$ is called an elementary
$r_j$-collapse.

\begin{defin}
A pure $d$-dimensional simplicial complex $K$ is called {\em
tree-like} or vertex collapsible if it is collapsible to a
$d$-simplex with the use of elementary $0$-collapses alone. In
other words $K$ is shellable and for each $j\geq 2$, the
intersection $F_j\cap(\bigcup_{i<j}~F_i)$ is a proper face of
$F_j$.
\end{defin}

In order to establish analogs of Propositions~\ref{prop:help} and
\ref{prop:transport} for complexes $Hom(K,L)$, we prove a result
which shows that elementary vertex collapsing  provides a good
substitute and a partial generalization for the concept of
``foldings'' of graphs used in \cite{BabsonKozlov1}
\cite{Kozl2004} in the theory of graph complexes $Hom(G,H)$.

\begin{prop}
\label{prop:vertex-collapse} Suppose that the simplicial complex
$K'$ is obtained from $K$ by an elementary vertex collapse. In
other words we assume that $K = \sigma\cup K',\, \sigma\cap
K'=\sigma'$, where $\sigma$ is a simplex in $K$ and $\sigma'$ a
facet of $\sigma$. Assume that $\sigma'$ is not maximal in $K'$,
i.e.\ that for some simplex $\sigma''\in K'$ and a vertex
$u\in\sigma'',\, \sigma' = \sigma''\setminus\{u\}$. Then for any
simplicial complex $L$, the inclusion map $\gamma : K'\rightarrow
K$ induces a homotopy equivalence
\[
\widehat{\gamma} : Hom(K,L)\rightarrow Hom(K',L).
\]
\end{prop}

\medskip\noindent
{\bf Proof:} Let $\{v\}=\sigma\setminus \sigma'$. Aside from the
inclusion map $\gamma : K'\rightarrow K$, there is a retraction
(folding) map $\rho : K\rightarrow K'$, where $\rho(v)=u$ and
$\rho\vert_{K'}=I_{K'}$. Since $\rho\circ\gamma = I_{K'}$, we
observe that $\widehat{\gamma}\circ\widehat{\rho}=Id_{K'}$ is the
identity map on $Hom(K',L)$, i.e.\ the complex $Hom(K',L)$ is a
retract of the complex $Hom(K,L)$. It remains to be shown that
$\widehat{\rho}\circ\widehat{\gamma}\simeq Id_{K}$ is homotopic to
the identity map on $Hom(K,L)$.

Note that if $\eta\in Hom(K,L)$ then $\eta':=
\widehat{\rho}\circ\widehat{\gamma}(\eta)$ is the function defined
by
 $$
 \eta'(w) = \left\{\begin{array}{rl} \eta(w), & \mbox{ {\rm
if} \, }
w\neq v \\
\eta(u), & \mbox{ {\rm if} \, } w = v.
\end{array}\right.
 $$
Let $\omega : Hom(K,L)\rightarrow Hom(K,L)$ be the map defined by
$$
 \omega(\eta)(w) = \left\{\begin{array}{ll} \eta(w), & \mbox{ {\rm
if} \, }
w\neq v \\
\eta(u)\cup\eta(v), & \mbox{ {\rm if} \, } w = v.
\end{array}\right.
 $$
Note that $\omega$ is well defined since if a vertex $x$ is
adjacent to $v$ it is also adjacent to $u$, hence the condition
$\omega(\eta)(v)\cap\eta(x)=\emptyset$ is a consequence of
$\eta(u)\cap\eta(x)=\emptyset = \eta(v)\cap\eta(x)$.

Since for each $\eta\in Hom(K,L)$ and each vertex $x\in K$,
\[
\eta(x)\subset \omega(\eta)(x) \supset
\widehat{\rho}\circ\widehat{\gamma}(\eta)(x),
\]
by the Order Homotopy Theorem \cite{Bjorner} \cite{Quillen}
\cite{Segal} all three maps $Id_{K}, \omega$ and
$\widehat{\rho}\circ\widehat{\gamma}$ are homotopic. This
completes the proof of the proposition. \hfill $\square$

\begin{cor}
If $T$ is a $d$-dimensional, tree-like simplicial complex than
$Hom(T,L)$ has the same homotopy type as the complex
$Hom(\Delta^d,L)$.
\end{cor}

\subsubsection{Parallel transport of homotopy types of maps}

As in the case of graph complexes, the real justification for the
introduction of the parallel transport of $Hom$-complexes comes
from the fact that it preserves the homotopy type of the maps
$Hom(K,L)\rightarrow Hom(\sigma, L)$. As a preliminary step we
prove an analogue of Proposition~\ref{prop:help}.

\begin{prop}\label{prop:helphelp}
Suppose that $\sigma_1$ and $\sigma_2$ are two distinct, adjacent
$k$-dimensional simplices in a finite simplicial complex $K$ which
share a common $(k-1)$-dimensional simplex $\tau$. Let $\Sigma =
\sigma_1\cup\sigma_2$. Let $\alpha : \Sigma\rightarrow\Sigma$ be
the automorphism of $\Sigma$ which interchanges simplices
$\sigma_1$ and $\sigma_2$ keeping the common face $\tau$
point-wise fixed.

Suppose that $\gamma_i: \sigma_i \rightarrow \Sigma$ is an obvious
embedding and $\widehat{\gamma}_i$ the associated maps of
$Hom$-complexes. Then,
 \begin{enumerate}
 \item[{\rm (a)}] the induced map $\widehat{\alpha}:
Hom(\Sigma,L)\rightarrow Hom(\Sigma,L)$ is homotopic to the
identity map $I_\Sigma$, and
 \item[{\rm (b)}]  the diagram
\[
\begin{CD}
Hom(\Sigma,L) @>=>> Hom(\Sigma,L) \\
 @V\widehat{\gamma}_{1}VV @VV\widehat{\gamma}_{2}V \\
 Hom(\sigma_1,L) @<<\mathcal{P}^H(\overrightarrow{\sigma_1\sigma_2})<
 Hom(\sigma_2,L)
 \end{CD}
\]
is commutative up to homotopy.
 \end{enumerate}
\end{prop}

\medskip\noindent
{\bf Proof:} By Proposition~\ref{prop:vertex-collapse}, both maps
$\widehat{\gamma}_i : Hom(\Sigma,L)\rightarrow Hom(\sigma_i,L)$
for $i=1,2$ are homotopy equivalences. Let $\rho_1 :
\Sigma\rightarrow \sigma_1$ and $\rho_2 : \Sigma\rightarrow
\sigma_2$ be the folding maps. Then $\rho_i\circ\gamma_i =
I_{\sigma_i},\, $ $\widehat{\gamma_i}\circ\widehat{\rho_i}= I$ and
we conclude that $\widehat{\rho}_i : Hom(\sigma_i,L)\rightarrow
Hom(\Sigma,L)$ is also a homotopy equivalence.

Part (a) of the proposition follows from the fact that
$\rho_1\circ\alpha\circ\gamma_1 = I_{\sigma_1}$ is an identity
map. Indeed, an immediate consequence is that
$\widehat{\gamma_1}\circ \widehat{\alpha}\circ \widehat{\rho_1} =
I : Hom(\sigma_1,L)\rightarrow Hom(\sigma_1,L)$ is also an
identity map and, in light of the fact that $\widehat{\gamma_1}$
and $\widehat{\rho_1}$ are homotopy inverses to each other, we
deduce that $\widehat{\alpha}\simeq I$.

For the part (b) we begin by an observation that $\rho_2\circ
\alpha\circ \gamma_1 = \overrightarrow{\sigma_1\sigma_2}$. Then,
$\mathcal{P}^H(\overrightarrow{\sigma_1\sigma_2}) =
\widehat{\gamma}_1\circ\widehat{\alpha}\circ\widehat{\rho}_2$, and
as a consequence of $\widehat{\alpha} \simeq I$ and the fact that
$\widehat{\rho}_2\circ \widehat{\gamma}_2\simeq I$, we conclude
that
\[
\mathcal{P}^H(\overrightarrow{\sigma_1\sigma_2})\circ
\widehat{\gamma}_2 = \widehat{\gamma}_1\circ\widehat{\alpha}\circ
\widehat{\rho}_2\circ\widehat{\gamma}_2\simeq \widehat{\gamma}_1.
\]
\hfill $\square$

\begin{prop}\label{prop:prenos}
Suppose that $K$ and $L$ are finite simplicial complexes and
$\sigma_1,\sigma_2$ a pair of adjacent (distinct), $k$-dimensional
simplices in $K$. Let $\alpha_i : \sigma_i\rightarrow K$ be the
embedding of $\sigma_i$ in $K$ and $\widehat{\alpha}_i:
Hom(K,L)\rightarrow Hom(\sigma_i,L)$ the associated map of
$Hom$-complexes. Then the following diagram commutes up to a
homotopy.

\begin{equation}\label{eqn:shema1}
\begin{CD}
Hom(K,L) @>=>> Hom(K,L) \\
 @V\widehat{\alpha}_{1}VV @VV\widehat{\alpha}_{2}V \\
 Hom(\sigma_1,L) @<<\mathcal{P}^L(\overrightarrow{\sigma_1\sigma_2})<
 Hom(\sigma_2,H)
 \end{CD}
\end{equation}
\end{prop}

\medskip\noindent
{\bf Proof:} Let $\Sigma := \sigma_1\cup\sigma_2, \, \tau :=
\sigma_1\cap\sigma_2$. Then $\alpha_i = \beta\circ\gamma_i$ where
$\beta:\Sigma\rightarrow K$ and $\gamma_i: \sigma_i \rightarrow
\Sigma$ are natural embeddings of complexes.

The diagram (\ref{eqn:shema1}) can be factored as
\begin{equation}\label{eqn:shema2}
\begin{CD}
  Hom(K,L) @>=>> Hom(K,L) \\
  @V\widehat{\beta} VV @VV\widehat{\beta} V \\
  Hom(\Sigma,L) @>=>> Hom(\Sigma,L)\\
  @V\widehat{\gamma}_1VV @VV\widehat{\gamma}_2 V\\
  Hom(\sigma_1,L) @<<\mathcal{P}^L(\overrightarrow{\sigma_1\sigma_2})<
  Hom(\sigma_2,L)
 \end{CD}
\end{equation}
Then the result is a direct consequence of
Proposition~\ref{prop:helphelp}, part (b).
 \hfill $\square$

\begin{cor}\label{cor:posledica}
Suppose that $K$ and $L$ are finite simplicial complexes, $\sigma$
a $k$-dimensional simplex in $K$ and $\alpha : \sigma\rightarrow
K$ the associated embedding. Let $\tau\in \Pi(K,\sigma)$. Then the
following diagram commutes up to a homotopy.

\begin{equation}\label{eqn:shema3}
\begin{CD}
Hom(K,L) @>=>> Hom(K,L) \\
 @V\widehat{\alpha} VV @VV\widehat{\alpha} V \\
 Hom(\sigma,L) @<<\widehat{\tau} <
 Hom(\sigma,L)
 \end{CD}
\end{equation}
\end{cor}

\subsection{A generalized B-K-L theorem} \label{sec:main}

In this section we prove the promised extension of the
Babson-Kozlov-Lov\'{a}sz theorem. The graphs are replaced by pure
$d$-dimensional simplicial complexes, while the role of the odd
cycle $C_{2r+1}$ is played by a complex $\Gamma$ which has some
special symmetry properties in the sense of the following
definition.

As usual, an involution $\omega : X\rightarrow X$ is the same as a
$\mathbb{Z}_2$-action on $X$. An involution on a simplicial
complex $\Gamma$ induces an involution on the complex
$Hom(\Gamma,L)$ for each simplicial complex $L$. For all other
standard facts and definitions related to $\mathbb{Z}_2$-complexes
the reader is referred to \cite{Matousek}.

\begin{defin}
\label{def:Phi-complex} A pure $d$-dimensional simplicial complex
$\Gamma$ is a $\Phi_d$-complex if it is a $\mathbb{Z}_2$-complex
with an invariant $d$-simplex $\sigma = \{v_0,v_1,\ldots, v_d\}$
such that the restriction $\tau :=\omega\vert_\sigma$ of the
involution $\omega : \Gamma\rightarrow\Gamma$ on $\sigma$ is a
non-trivial element of the group $\Pi(\Gamma,\sigma)$.
\end{defin}

\begin{rem}{\rm
By definition, if $\Gamma$ is a $\Phi_d$-complex then the
inclusion map $\alpha: \sigma\rightarrow\Gamma$ is
$\mathbb{Z}_2$-equivariant, so the associated map
$\widehat{\alpha}:Hom(\Gamma,K)\rightarrow Hom(\sigma,K)$ is also
$\mathbb{Z}_2$-equivariant for each complex $K$. }
\end{rem}

\begin{figure}[hbt]
\centering
\includegraphics[scale=0.40]{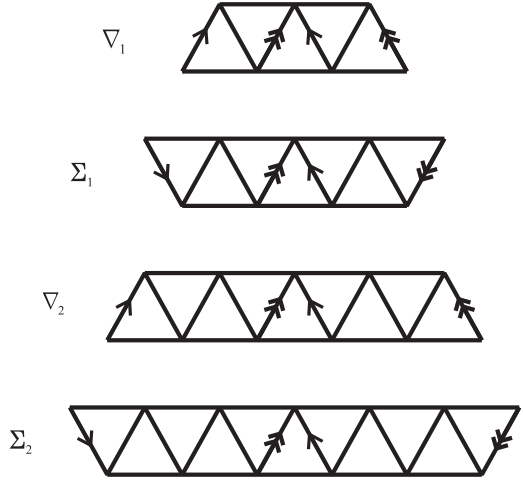}
\caption{Examples of $\Phi_2$-complexes.} \label{fig:piramida}
\end{figure}

\begin{exam}{\rm
The graph $C_{2r+1}$ is obviously an example of a
$\Phi_1$-complex. Figure~\ref{fig:piramida} displays four examples
of $\Phi_2$-complexes, initial elements of two infinite series
$\nabla_\mu$ and $\Sigma_\nu,\, \mu,\nu\in \mathbb{N}$. The
complexes $\nabla_1$ and $\nabla_2$ etc.\ are obtained from two
triangulated annuli, glued together along a common triangle
$\sigma$. Similarly, the complexes $\Sigma_1, \Sigma_2, \ldots$,
are obtained by gluing together two triangulated M\"{o}bius
strips. The associated group of projectivities are
$\Pi(\nabla_\mu,\sigma)=S_3$ and
$\Pi(\Sigma_\nu,\sigma)=\mathbb{Z}_2$. }
\end{exam}

\begin{theo}\label{thm:main}
Suppose that $\Gamma$ is a $\Phi_d$-complex in the sense of
Definition~\ref{def:Phi-complex}, with an associated invariant
simplex $\sigma = \{v_0,v_1,\ldots, v_d\}$. Suppose that $K$ is a
pure $d$-dimensional simplicial complex. Than for $m$ even,
\begin{equation}\label{eq:complex}
Coind_{\mathbb{Z}_2}(Hom(\Gamma,K))\geq m \Rightarrow \chi(K)\geq
m+d+2.
\end{equation}
\end{theo}

\medskip\noindent
{\bf Proof:} By definition
$Coind_{\mathbb{Z}_2}(Hom(\Gamma,K))\geq m$ means that there
exists a $\mathbb{Z}_2$-equiva\-ri\-ant map $\mu: S^m\rightarrow
Hom(\Gamma, K)$. Assume that $\chi(K)\leq m+d+1$ which means that
there exists a non-degenerate simplicial map $\phi :
K\rightarrow\Delta^{[m+d+1]}$. By functoriality of the
construction of $Hom$-complexes, Section~\ref{sec:naturality},
there is an induced $\mathbb{Z}_2$-equivariant map
$\widehat{\phi}: Hom(\Gamma, K)\rightarrow Hom(\Gamma,
\Delta^{[m+d+1]})$ and similarly a map $\widehat{\alpha}:
Hom(\Gamma, \Delta^{[m+d+1]})\rightarrow Hom(\sigma,
\Delta^{[m+d+1]})$. By \cite{Kozl2004} Theorem~3.3.3., the complex
 \[
Hom(\sigma, \Delta^{[m+d+1]})\cong Hom(K_{d+1},K_{m+d+1})
 \]
is a wedge of $m$-dimensional spheres. Since $Hom(\sigma,
\Delta^{[m+d+1]})$ is a free $\mathbb{Z}_2$-complex, we deduce
that there exists a $\mathbb{Z}_2$-map $Hom(\sigma,
\Delta^{[m+d+1]})\rightarrow S^m$. All these maps can be arranged
in the following sequence of $\mathbb{Z}_2$-equivariant maps
 \[
S^m   \stackrel{\mu}\longrightarrow Hom(\Gamma, K)
\stackrel{\widehat{\phi}}\longrightarrow Hom(\Gamma,
\Delta^{[m+d+1]}) \stackrel{\widehat{\alpha}}\longrightarrow
Hom(\sigma, \Delta^{[m+d+1]}) \stackrel{\nu}\longrightarrow S^m.
 \]
By Corollary~\ref{cor:posledica}, there is a homotopy equivalence
$\widehat{\alpha}\simeq \tau\circ \widehat{\alpha}$. This
contradicts Proposition~\ref{prop:not-homotopic}, which completes
the proof of the theorem. \hfill $\square$

\begin{prop}\label{prop:not-homotopic}
Suppose that $f: X\rightarrow Y$ is a $\mathbb{Z}_2$-equivariant
map of free $\mathbb{Z}_2$-complexes $X$ and $Y$ where
$\mathbb{Z}_2=\{1,\omega\}$. Assume that
$Coind_{\mathbb{Z}_2}(X)\geq m\geq Ind_{\mathbb{Z}_2}(Y)$, where
$m$ is an even integer. In other words our assumption is that
there exist $\mathbb{Z}_2$-equivariant maps $\mu$ and $\nu$ such
that
\[
S^m \stackrel{\mu}\longrightarrow X \stackrel{f}\longrightarrow Y
\stackrel{\nu}\longrightarrow S^m .
\]
Then the maps $f$ and $\omega\circ f$ are not homotopic.

\medskip\noindent
{\bf Proof:} If $f\simeq \omega\circ f : X\rightarrow Y$ then
$\nu\circ f\circ\mu \simeq \nu\circ\omega\circ f\circ\mu :
S^m\rightarrow S^m$ and by the equivariance of $\nu\,$,
$\omega\circ g\simeq g : S^m\rightarrow S^m$ where $g:= \nu\circ
f\circ \mu$. It follows that
\[
-{\rm deg}(g) = {\rm deg}(\omega) {\rm deg}(g) = {\rm
deg}(\omega\circ g) = {\rm deg}(g),
\]
i.e.\ ${\rm deg}(g)=0$, which is in contradiction with a well
known fact \cite{Matousek} that a $\mathbb{Z}_2$-equivariant map
$g : S^m\rightarrow S^m$ of spheres must have an odd degree.
\hfill $\square$
\end{prop}

\begin{cor}\label{cor:LBK-gen}
Suppose that $\Gamma$ is a $\Phi_d$-complex with an associated
invariant simplex $\sigma = \{v_0,v_1,\ldots, v_d\}$. Suppose that
$K$ is a pure $d$-dimensional simplicial complex. Than for $k$
odd,
\begin{equation}\label{eq:za-kraj}
Hom(\Gamma,K) \mbox{ {\rm is $k$-connected} } \Rightarrow
\chi(K)\geq k+d+3.
\end{equation}
\end{cor}

\medskip\noindent
{\bf Proof:} If $Hom(\Gamma,K)$ is $k$-connected then $Coind_{
\mathbb{Z}_2}(Hom(\Gamma,K))\geq k+1$, hence the implication
(\ref{eq:za-kraj}) is an immediate consequence of
Theorem~\ref{thm:main}. \hfill $\square$

\subsection{Flag complexes and the homotopy test graphs} \label{sec:paradighm}

A simplicial complex is a flag-complex if $\sigma\in K$ if and
only if the associated $1$-skeleton $\sigma^{(1)}$ of $\sigma$ is
a subcomplex of $K$. Examples of flag-complexes include the
clique-complex $Clique(G)$ of a graph and the order complex
$\Delta(P)$ of all chains (flags) in a poset $P$.

If $K$ is a flag-complex and $G=G_K:= K^{(1)}$ is the associated
vertex-edge graph of $K$, then $K = Clique(G)$. From here we
deduce, relying on the isomorphism given in
Example~\ref{exam:companion}, that if $K$ and $L$ are flag
complexes then

 \begin{equation}\label{eq:useful}
Hom(K,L) \cong Hom(G_K, G_L).
 \end{equation}
This simple observation allows us to transfer results about
simplicial complexes back to graphs. An example is the following
theorem which can be deduced from Corollary~\ref{cor:LBK-gen}.
\begin{theo}\label{thm:paradigm} Suppose that $W$ is $\Phi_d$-complex
(Definition~\ref{def:Phi-complex}) which is also a flag complex.
Assume that $k$ is an odd integer.  Let $T:=W^{(1)}$. Then for
each graph $G$
\begin{equation}\label{eq:paradigm}
Hom(T,G) \mbox{ {\rm is $k$-connected} } \Rightarrow \chi(G)\geq
k+d+3.
\end{equation}
\end{theo}

Examples of  $\Phi_2$-complexes $W$ satisfying conditions of
Theorem~\ref{thm:paradigm} are the complexes
$\Sigma_1,\nabla_2,\Sigma_2$ depicted in
Figure~\ref{fig:piramida}. Moreover, for all associated graphs
$\chi(W^{(1)})=4$ which in light of (\ref{eq:paradigm}) means that
they are all candidates for ``homotopy test graphs'' in the sense
of \cite[Chapter~6]{Kozlov-review}.

\end{document}